\documentclass[letterpaper,conference]{ieeeconf}

\IEEEoverridecommandlockouts

\usepackage{amsmath}
\usepackage{amssymb}
\usepackage{mathrsfs}
\usepackage{graphicx}
\usepackage{amssymb}
\usepackage{esvect}
\usepackage{bbm}
\usepackage{algorithmic}
\usepackage{algorithm}
\usepackage{url}

\newtheorem{lemma}{Lemma}
\newtheorem{remark}{Remark}

\DeclareMathOperator{\step}{step}

\usepackage{bbm,bm}

\newcommand{\longversion}{}

\title{Scheduling~Rigid~Demands~on~Continuous-Time~Linear~Shift-Invariant~Systems \thanks{This work was supported by the Australian Research Council (LP130100605), Rubicon Water~Pty~Ltd, and a McKenzie Fellowship.}}
\author{Farhad Farokhi, Michael Cantoni, and Iman Shames\thanks{The authors are with the Department of Electrical and Electronic Engineering, The University of Melbourne, Parkville, Victoria 3010, Australia. Emails:\{ffarokhi,cantoni,ishames\}@unimelb.edu.au}}

\begin{document}
\maketitle

\begin{abstract} We consider load scheduling on constrained continuous-time linear dynamical systems, such as automated irrigation and other distribution networks. The requested loads are  rigid, i.e.,  the shapes cannot be changed. Hence, it is only possible to shift the order back-and-forth in time to arrive at a feasible schedule. We present a numerical algorithm based on using log-barrier functions to include the state constraints in the social cost function (i.e., an appropriate function of the scheduling delays). This algorithm requires a feasible initialization. Further, in another algorithm, we treat the state constraints as soft constraints and heavily penalize the constraint violations. This algorithm can even be initialized at an infeasible point. The applicability of both these numerical algorithms is demonstrated on an automated irrigation network with two pools and six farms. 
\end{abstract}

\section{Introduction}
Scheduling problems arise in a variety of contexts. A peculiar scheduling problem is studied in this paper. It involves a constrained dynamical system and the processing of request to apply load, with a fixed but shiftable profile, on this system across time. The goal is to optimize a social measure of sensitivity to scheduling delay while satisfying hard constraints. This problem is motivated by an aspect of demand management in automated irrigation networks~\cite{5400193,hong2012optimization, reddy1999optimal}, and may arise in other areas. The main challenge associated with this problem relates to the rigidity of the load request, whereby the construction of a feasible schedule can only involve shifting requests back-and-forth in time. Relaxation of the rigidity requirement can lead to a formulation as a (large) linear program~\cite{5400193}.

The formulation of the rigid load scheduling problem here distinguishes itself in the following ways. By contrast with~\cite{hong2012optimization, reddy1999optimal}, a dynamics relationship between the load and the constrained system states are modelled. In~\cite{hong2012optimization, reddy1999optimal}, only static capacity constraints are considered. The load scheduling problem considered in~\cite{5400193} does include dynamics, however this is modelled in discrete time. By contrast a continuous-time setting is employed in this paper. The discrete time formulation in~\cite{5400193} leads to a mixed-integer program, which is difficult to solve~\cite{garey1979computers}. The continuous-time formulation here, on the other hand, gives rise to two gradient based numerical algorithms. The first involve log-barrier functions and thus a feasible initial point. The other uses a soft encoding of the state constraints, with heavy penalty on constraint violation, which does not require a feasible initial point. Both algorithms lead to only locally optimal solutions due to the non-convexity of the scheduling problem.

The rest of the paper is organized as follows. We first formulate the problem in Section~\ref{sec:problem}. The numerical algorithms are presented in Sections~\ref{sec:numerical} and Section~\ref{sec:soft}. In Section~\ref{sec:numerical}, the applicability of the developed algorithms is numerically studied on an automated irrigation network with two pools and six farms. Finally, Section~\ref{sec:conclusions} concludes the paper.

\section{Problem Formulation} \label{sec:problem}
Consider the continuous-time linear time-invariant dynamical system 
\begin{align} \label{eqn:1}
\dot{x}(t)=Ax(t)+Bu(t)+\sum_{i=1}^m E_iw_i(t),\,\,x(0)=x_0,
\end{align}
where $x(t)\in\mathbb{R}^{n_x}$ is the state of the system, $u(t)\in\mathbb{R}^{n_u}$ is the control input (e.g., the water-level references in automated irrigation networks), and $w_i(t)\in\mathbb{R}^{n_{w,i}}$, $1\leq i\leq m$, is a profile-constrained (e.g. on-off) input signal representing the scheduled load on the system corresponding to the supply of resources to customer $i$ (e.g., the flow of the supplied water to each farmer in an irrigation network). Throughout this paper, we assume that the control signal $u(t)$ over the planning horizon $[0,T]$ with $T\in\mathbb{R}_{>0}$ is a discrete-time signal passed through a zero-order hold, that is, $u_i(t)=\alpha_{i,k}$ for all $1\leq i\leq n_u$ and all $k\Delta\leq t< (k+1)\Delta$ with a given sampling time $0<\Delta\leq T$. Although slightly conservative, this assumption allows us to work with finite-dimensional optimization problems instead of more complicated optimal control problems. For all integers $0\leq k\leq K:=\lceil T/\Delta\rceil-1$ and $1\leq i\leq n_u$, we define 
$\xi_{i,k}(t)=e_i[\step(t-k\Delta)-\step(t-(k+1)\Delta)],\forall t\in[0,T],$
where the mapping $\step:\mathbb{R}\rightarrow \{0,1\}$ denotes the Heaviside step function, i.e., $\step(t)=1$ if $t\geq 0$ and $\step(t)=0$ otherwise. Moreover, $e_i\in\mathbb{R}^{n_u}$ is the column-vector with all entries equal to zero except the $i$-th entry which is equal to one. Therefore, we get
$u(t)=u_0+\sum_{k=1}^K\sum_{i=1}^{n_u} \alpha_{i,k}\xi_{i,k}(t), \forall t\in[0,T],$
where $u_0\in\mathbb{R}^{n_u}$ is the steady-state control input.

The customers submit demands $(v_i(t))_{t\in\mathbb{R}}$, $1\leq i\leq m$. These demands are  rigid (i.e., their shape cannot be changed). Hence, our decision variables are the delays that correspond to shifting the requested demand across the planning horizon, i.e., we  select $\tau_i>0$ so that $w_i(t)=v_i(t-\tau_i)$ for each $1\leq i\leq m$. In doing so, the goal is to ensure that that the state of the network $x(t)$ stays inside the feasible set $\mathcal{X}=\{x\in\mathbb{R}^{n_x}\,|\,Cx\leq d\}$. We can write this scheduling problem as
\begin{subequations}\label{eqn:optim:1}
\begin{align} 
\min_{(\tau_i)_{i=1}^m,((\alpha_i)_{i=1}^{n_u})_{k=1}^K} \, & \sum_{i=1}^m h_i(\tau_i),\\[-.5em]
\mathrm{s.t.} \hspace{.40in} & \underline{\tau}_i\leq \tau_i\leq  \overline{\tau}_i, \,\forall i\in\{1,\dots,m\} \\[-.2em]
& \dot{x}(t)=Ax(t)+Bu(t)\nonumber\\[-.5em]
&\hspace{.16in}+\sum_{i=1}^m E_iv_i(t-\tau_i),x(0)=x_0,\\[-.4em] 
& Cx(t)\leq d, \forall t\in[0,T], \label{eqn:optim:1:sc}\\[-.5em]
& u(t)\hspace{-.04in}=\hspace{-.04in}\sum_{k=1}^K\hspace{-.04in}\sum_{i=1}^{n_u} \hspace{-.04in}\alpha_{i,k}\xi_{i,k}(t), \forall t\in[0,T],\\[-.3em]
& u_0\hspace{-.03in}+\hspace{-.03in}\underline{u}\hspace{-.03in}\leq\hspace{-.03in} u(t)\hspace{-.03in}\leq\hspace{-.03in} u_0\hspace{-.03in}+\hspace{-.03in}\overline{u}, \forall t\in[0,T],
\end{align}
\end{subequations}
where $\underline{\tau}_i$ and $\overline{\tau}_i$ are the bounds on the scheduling delay for demand $i$, $\underline{u}$ and $\overline{u}$ are the bounds on the control signal deviations $u(t)-u_0$, and the continuously differentiable mapping $h_i:\mathbb{R}\rightarrow\mathbb{R}$ captures the sensitivity of customer $i$ to the delay for scheduling its demand. Throughout the next section, we implicitly assume that $T$ is long enough so that the optimization problem in~\eqref{eqn:optim:1} becomes feasible with a constant nominal control input (i.e., if the demands are separated from each other ``to some degree'', the state of the system stays feasible without any effort). This assumption is made to make sure that we can always find a feasible initial condition for the numerical algorithm, proposed in the next section, by simply separating the demands from each other. Towards the end of this paper, we present another approach for solving our scheduling problem that avoids requiring a feasible initial condition by treating the constraints on the state as soft constraints.

\section{Numerical Algorithm} \label{sec:algorithm}
In this section, we present a numerical algorithm for solving~\eqref{eqn:optim:1} by adding the state constraints in~\eqref{eqn:optim:1:sc} to the cost function using log-barrier functions. Let us define
\begin{align*}
\bar{x}_0(t)&=\exp(At)x_0+\int_0^t \exp(A(t-\beta))Bu_0\mathrm{d}\beta,\\[-.4em]
\bar{x}_{i,k}^u(t)&=\int_0^t \exp(A(t-\beta))B\xi_{i,k}(\beta)\mathrm{d}\beta, \forall i\in\{1,\dots,n_u\},\\[-.8em]
&\hspace{1.97in}\forall k\in\{1,\dots,K\},\\[-.8em]
\bar{x}_i^v(t)&=\int_0^t \exp(A(t-\beta))E_iv_i(\beta)\mathrm{d}\beta, \forall i\in\{1,\dots,m\}.
\end{align*}
Since the underlying system in~\eqref{eqn:1} is linear and time invariant, the solution of the ordinary differential equation~\eqref{eqn:1} can be written explicitly as
$
x(t)=\bar{x}_0(t)+\sum_{k=1}^K\sum_{i=1}^{n_u}\alpha_{i,k}\bar{x}^u_{i,k}(t)+\sum_{i=1}^m\bar{x}_i^v(t-\tau_i).
$
Now, we can rewrite the optimization problem in~\eqref{eqn:optim:1} as
\begin{subequations}\label{eqn:optim:2}
\begin{align} 
\min_{(\tau_i)_{i=1}^m,((\alpha_{i,k})_{i=1}^{n_u})_{k=1}^K} \, & \sum_{i=1}^m h_i(\tau_i),\\[-1em]
\mathrm{s.t.} \hspace{.44in} & x(t)=\bar{x}_0(t)+\sum_{k=1}^K\sum_{i=1}^{n_u}\alpha_{i,k}\bar{x}^u_{i,k}(t)\nonumber\\[-.4em]
&\hspace{.72in}+\sum_{i=1}^m\bar{x}_i^v(t-\tau_i),\\[-.3em]
& Cx(t)\leq d, \forall t\in[0,T], \label{eqn:optim:2:sc}\\[-.3em]
& \underline{\tau}_i\leq \tau_i\leq  \overline{\tau}_i, \,\forall i\in\{1,\dots,m\},\\[-.3em]
& \underline{u}_i \leq  \alpha_{i,k}\leq \overline{u}_i, \forall i\in\{1,\dots,n_u\},\nonumber\\[-.3em]
&\hspace{0.91in}\forall k\in\{1,\dots,K\}.
\end{align}
\end{subequations}
This optimization problem is still difficult to solve as we have to check infinitely many constraints; see~\eqref{eqn:optim:2:sc}. Let us use the notation $C_j$, $1\leq j\leq p$, to denote the rows of the matrix $C\in\mathbb{R}^{p\times n_x}$. We add the state constraints in~\eqref{eqn:optim:2:sc} to the cost function using log-barrier functions. This transforms the optimization problem in~\eqref{eqn:optim:2} to
\begin{subequations}\label{eqn:optim:4}
\begin{align} 
\min_{(\tau_i)_{i=1}^m,((\alpha_{i,k})_{i=1}^{n_u})_{k=1}^K} & J((\tau_i)_{i=1}^m,((\alpha_{i,k})_{i=1}^{n_u})_{k=1}^K)\\
\mathrm{s.t.} \hspace{.43in} & \underline{\tau}_i\leq \tau_i\leq  \overline{\tau}_i, \,\forall i\in\{1,\dots,m\},\\
& \underline{u}_i \leq  \alpha_{i,k}\leq \overline{u}_i, \forall i\in\{1,\dots,n_u\},\nonumber\\
&\hspace{0.91in}\forall k\in\{1,\dots,K\},
\end{align}
\end{subequations}
where 
\begin{align*}
J((\tau_i&)_{i=1}^m,((\alpha_{i,k})_{i=1}^{n_u})_{k=1}^K)\\[-.3em]
=&\sum_{i=1}^m h_i(\tau_i)-\sum_{z=1}^p\int_0^T \epsilon\log\bigg(-C_z\bigg[\bar{x}_0(t)\\[-.5em]&+\sum_{k=1}^K\sum_{i=1}^{n_u}\alpha_{i,k}\bar{x}^u_{i,k}(t)+\sum_{i=1}^m\bar{x}_i^v(t-\tau_i)\bigg]+d_z\bigg) \mathrm{d}t
\end{align*} 
in which $\epsilon\in\mathbb{R}_{>0}$ is an appropriately selected parameter. 

\begin{remark} With increasing $\epsilon$, the optimal solution is pushed further from the boundary of the feasible set. Therefore, to recover the optimal scheduling, we need to sequentially reduce $\epsilon$ and employ the solution of each step as the initialization of the next step. This would result in a  more numerically stable algorithm; see the log barrier methods in~\cite{boyd2004convex}.
\end{remark}

\begin{lemma} \label{lemma:1} $J((\tau_i)_{i=1}^m,((\alpha_{i,k})_{i=1}^{n_u})_{k=1}^K)$ is a continuously differentiable function. Moreover, 
\begin{align*}
\frac{\partial}{\partial \tau_\ell} & J((\tau_i)_{i=1}^m,((\alpha_{i,k})_{i=1}^{n_u})_{k=1}^K)=\frac{\mathrm{d}}{\mathrm{d} \tau_\ell} h_\ell(\tau_\ell)\\[-.3em]
+&\sum_{z=1}^p\int_0^T \epsilon\frac{-C_z(A\bar{x}_\ell^v(t-\tau_\ell)+E_\ell v_\ell(t-\tau_\ell))}{-C_z x(t;(\tau_i)_{i=1}^m,((\alpha_{i,k})_{i=1}^{n_u})_{k=1}^K)+d_z}\mathrm{d}t\\[-.3em]
\frac{\partial}{\partial \alpha_{j,\ell}}& J((\tau_i)_{i=1}^m,((\alpha_{i,k})_{i=1}^{n_u})_{k=1}^K)\\[-.3em]
=&\sum_{z=1}^p\int_0^T \epsilon\frac{C_z\bar{x}_{j,\ell}^u(t)}{-C_zx(t;(\tau_i)_{i=1}^m,((\alpha_{i,k})_{i=1}^{n_u})_{k=1}^K)+d_z}\mathrm{d}t
\end{align*}
where
\begin{align*}
x(t;(&\tau_i)_{i=1}^m,((\alpha_{i,k})_{i=1}^{n_u})_{k=1}^K)\\[-.3em]
&=\bar{x}_0(t)+\sum_{k=1}^K\sum_{i=1}^{n_u}\alpha_{i,k}\bar{x}^u_{i,k}(t)+\sum_{i=1}^m\bar{x}_i^v(t-\tau_i).
\end{align*}
\end{lemma}

\ifdefined\shortversion {\color{blue}
\begin{proof}
See~\cite{reportwithproof} for a detailed proof.  
\end{proof} }
\fi
\ifdefined\longversion
\begin{proof} First note that
\begin{align*}
\frac{\partial}{\partial \tau_\ell}&\int_0^T \epsilon\log(-C_z x(t;(\tau_i)_{i=1}^m,((\alpha_{i,k})_{i=1}^{n_u})_{k=1}^K)+d_z)\mathrm{d}t\\
&=\int_0^T \epsilon\frac{-C_z\partial x(t;(\tau_i)_{i=1}^m,((\alpha_{i,k})_{i=1}^{n_u})_{k=1}^K)/\partial \tau_\ell}{-C_z x(t;(\tau_i)_{i=1}^m,((\alpha_{i,k})_{i=1}^{n_u})_{k=1}^K)+d_z}\mathrm{d}t
\end{align*}
where
\begin{align*}
\frac{\partial}{\partial \tau_\ell}x(t;(\tau_i)_{i=1}^m,&((\alpha_{i,k})_{i=1}^{n_u})_{k=1}^K)\\
&=\frac{\partial}{\partial \tau_\ell}\bar{x}_\ell^v(t-\tau_\ell)\\
&=-\dot{\bar{x}}_\ell^v(t-\tau_\ell)\\
&=-(A\bar{x}_\ell^v(t-\tau_\ell)+E_\ell v_\ell(t-\tau_\ell)).
\end{align*}
Similarly, we have
\begin{align*}
\frac{\partial}{\partial \alpha_{j,\ell}}&\int_0^T \epsilon\log(-C_z x(t;(\tau_i)_{i=1}^m,((\alpha_{i,k})_{i=1}^{n_u})_{k=1}^K)+d_z)\mathrm{d}t\\
&= \int_0^T \epsilon\frac{-C_z\partial x(t;(\tau_i)_{i=1}^m,((\alpha_{i,k})_{i=1}^{n_u})_{k=1}^K)/\partial \alpha_{j,\ell}}{-C_z x(t;(\tau_i)_{i=1}^m,((\alpha_{i,k})_{i=1}^{n_u})_{k=1}^K)+d_z}\mathrm{d}t
\end{align*}
where
\begin{align*}
\frac{\partial}{\partial \alpha_{j,\ell}}x(t;(\tau_i)_{i=1}^m,((\alpha_{i,k})_{i=1}^{n_u})_{k=1}^K)
&=\bar{x}_{j,\ell}^u(t).
\end{align*}
The rest of the proof follows from simple algebraic manipulations.
\end{proof}
\fi

\begin{algorithm*}[t]
\caption{\label{alg:1} Projected gradient algorithm for scheduling rigid demands.}
\begin{algorithmic}[1]
\REQUIRE Feasible initialization $(\tau_i[0])_{i=1}^m$ and $((\alpha_{i,k}[0])_{i=1}^{n_u})_{k=1}^K$
\FOR{$l=1,2,\dots$}
\STATE{ Update
\begin{align*}
\hspace{-.6in}
\tau_\ell[l]=
P_{\underline{\tau}_\ell}^{\overline{\tau}_\ell}\left[\tau_\ell[l-1]
-\mu_l^{\tau_i} \frac{\partial J((\tau_i)_{i=1}^m,((\alpha_{i,k})_{i=1}^{n_u})_{k=1}^K))}{\partial \tau_i}  \bigg|_{\scriptsize
\begin{array}{c}
(\tau_i)_{i=1}^m=(\tau_i[l-1])_{i=1}^m\\
((\alpha_{i,k})_{i=1}^{n_u})_{k=1}^K=((\alpha_{i,k}[l-1])_{i=1}^{n_u})_{k=1}^K
\end{array}
}\right],\,\forall \ell\in\{1,\dots,m\},
\end{align*}
and
\begin{align*}
\hspace{-.6in}
\alpha_{j,\ell}[l]=&P_{\underline{u}_j}^{\overline{u}_j}\left[
\alpha_{j,\ell}[l-1]-\mu_l^{\alpha_{j,\ell}} \frac{\partial J((\tau_i)_{i=1}^m,((\alpha_{i,k})_{i=1}^{n_u})_{k=1}^K))}{\partial \alpha_{j,\ell}}  \bigg|_{\scriptsize
\begin{array}{c}
(\tau_i)_{i=1}^m=(\tau_i[l-1])_{i=1}^m\\
((\alpha_{i,k})_{i=1}^{n_u})_{k=1}^K=((\alpha_{i,k}[l-1])_{i=1}^{n_u})_{k=1}^K
\end{array}
}\right], \forall j\in\{1,\dots,n_u\},\nonumber\\
&\hspace{5.28in}\forall \ell\in\{1,\dots,K\},
\end{align*}
where, for constants $\beta<\gamma$, $P_\beta^\gamma[x]=\beta$ if $x<\beta$, $P_\beta^\gamma[x]=x$ if $\beta\leq x\leq \gamma$, and $P_\beta^\gamma[x]=\gamma$ if $x>\gamma$.
}
\ENDFOR
\end{algorithmic}
\end{algorithm*}

Now, we can use Algorithm~\ref{alg:1} (overleaf) to recover a local solution of~\eqref{eqn:optim:4}. We can select the step sizes $\mu_l^{\tau_i}$ and $\mu_l^{\alpha_{j,\ell}}$ using backtracking line search algorithm~\cite[p.\,464]{boyd2004convex} and terminate the algorithm whenever the improvements in the cost function becomes negligible. Unfortunately, this algorithm requires a feasible starting point (because the argument of the logarithmic functions cannot become negative). We remove this assumption in the next section by proposing a numerical procedure that treats the state constraints as soft constraints.

\section{Soft Constraints on States} \label{sec:soft}
In the previous section, we were required to find a feasible initialization to be able to run Algorithm~\ref{alg:1}. Here, we take a different approach by solving the optimization problem 
\begin{subequations}\label{eqn:optim:5}
\begin{align} 
\min_{(\tau_i)_{i=1}^m,((\alpha_i)_{i=1}^{n_u})_{k=1}^K} \, & \sum_{i=1}^m h_i(\tau_i)+\sum_{z=1}^p \int_0^Te^{\vartheta(C_zx(t)-d_z)}\mathrm{d}t,\\[-.3em]
\mathrm{s.t.} \hspace{.42in} & \underline{\tau}_i\leq \tau_i\leq  \overline{\tau}_i, \,\forall i\in\{1,\dots,m\} \\[-.3em]
& \dot{x}(t)=Ax(t)+Bu(t)\nonumber\\[-.3em]
&\hspace{.76in}+\sum_{i=1}^m E_iv_i(t-\tau_i), \\[-.3em]
&x(0)=x_0,\\[-.3em]
& u(t)\hspace{-.04in}=\hspace{-.04in}\sum_{k=1}^K\hspace{-.04in}\sum_{i=1}^{n_u}\hspace{-.04in} \alpha_{i,k}\xi_{i,k}(t), \forall t\in[0,T],\\[-.3em]
& \underline{u}\leq u(t)\leq \overline{u}, \forall t\in[0,T],
\end{align}
\end{subequations}
where $\vartheta\in\mathbb{R}_{>0}$ is an appropriately selected constant. In this problem, we may violate the constraints $Cx(t)-d\leq 0$, however, the term $\sum_{z=1}^p \int_0^Te^{\vartheta(C_zx(t)-d_z)}\mathrm{d}t$ heavily penalizes such violations. For small values of $\vartheta$, this term also penalizes the states being close to the boundary (of the feasible set), however, as we increase $\vartheta$, this term approaches zero inside the feasible set and infinity outside of the feasible set. 

Note that similar to the previous section, we can transform~\eqref{eqn:optim:5} into 
\begin{subequations}\label{eqn:optim:6}
\begin{align} 
\min_{(\tau_i)_{i=1}^m,((\alpha_{i,k})_{i=1}^{n_u})_{k=1}^K} \, & \sum_{i=1}^m h_i(\tau_i)+\sum_{z=1}^p \int_0^Te^{\vartheta(C_zx(t)-d_z)}\mathrm{d}t,\\[-.3em]
\mathrm{s.t.} \hspace{.40in} & x(t)=\bar{x}_0(t)+\sum_{k=1}^K\sum_{i=1}^{n_u}\alpha_{i,k}\bar{x}^u_{i,k}(t)\nonumber\\[-.3em]
&\hspace{.72in}+\sum_{i=1}^m\bar{x}_i^v(t-\tau_i),\\[-.3em]
& \underline{u}_i \leq  \alpha_{i,k}\leq \overline{u}_i, \forall i\in\{1,\dots,n_u\},\nonumber\\[-.3em]
&\hspace{0.91in}\forall k\in\{1,\dots,K\}.
\end{align}
\end{subequations}
Let us define
\begin{align*}
J'((\tau_i&)_{i=1}^m,((\alpha_{i,k})_{i=1}^{n_u})_{k=1}^K)\\
=&\sum_{i=1}^m h_i(\tau_i)+\sum_{z=1}^p \int_0^T\exp\bigg(\vartheta\bigg(C_z\bigg[\bar{x}_0(t)\\&+\sum_{k=1}^K\sum_{i=1}^{n_u}\alpha_{i,k}\bar{x}^u_{i,k}(t)+\sum_{i=1}^m\bar{x}_i^v(t-\tau_i)\bigg]-d_z\bigg)\bigg)\mathrm{d}t.
\end{align*} 
Hence, we may rewrite~\eqref{eqn:optim:6} as
\begin{subequations}\label{eqn:optim:7}
\begin{align} 
\min_{(\tau_i)_{i=1}^m,((\alpha_{i,k})_{i=1}^{n_u})_{k=1}^K} \, & J'((\tau_i)_{i=1}^m,((\alpha_{i,k})_{i=1}^{n_u})_{k=1}^K),\\
\mathrm{s.t.} \hspace{.45in} & \underline{\tau}_i\leq \tau_i\leq  \overline{\tau}_i, \,\forall i\in\{1,\dots,m\},\\ & \underline{u}_i \leq  \alpha_{i,k}\leq \overline{u}_i, \forall i\in\{1,\dots,n_u\},\nonumber\\
&\hspace{0.91in}\forall k\in\{1,\dots,K\}.
\end{align}
\end{subequations}
Similarly, we can prove the following result regarding the augmented cost function.

\begin{lemma} \label{lemma:2} $J'((\tau_i)_{i=1}^m,((\alpha_{i,k})_{i=1}^{n_u})_{k=1}^K)$ is a continuously differentiable function. Moreover,
\begin{align*}
&\frac{\partial}{\partial \tau_\ell}  J'((\tau_i)_{i=1}^m,((\alpha_{i,k})_{i=1}^{n_u})_{k=1}^K)=\frac{\mathrm{d}}{\mathrm{d} \tau_\ell} h_\ell(\tau_\ell)\\
&-\sum_{z=1}^p\int_0^T\vartheta\exp(\vartheta(C_zx(t;(\tau_i)_{i=1}^m,((\alpha_{i,k})_{i=1}^{n_u})_{k=1}^K)-d_z))\\
&\hspace{.4in}\times C_z(A\bar{x}_\ell^v(t-\tau_\ell)+E_\ell v_\ell(t-\tau_\ell))\mathrm{d}t,\\
&\frac{\partial}{\partial \alpha_{j,\ell}} J'((\tau_i)_{i=1}^m,((\alpha_{i,k})_{i=1}^{n_u})_{k=1}^K)\\
=&\sum_{z=1}^p\int_0^T\vartheta\exp(\vartheta(C_zx(t;(\tau_i)_{i=1}^m,((\alpha_{i,k})_{i=1}^{n_u})_{k=1}^K)-d_z))\\
&\hspace{.4in}\times C_z \bar{x}_{j,\ell}^u(t)\mathrm{d}t,
\end{align*}
where $x(t;(\tau_i)_{i=1}^m,((\alpha_{i,k})_{i=1}^{n_u})_{k=1}^K)$ is defined as in Lemma~\ref{lemma:1}.
\end{lemma}

\ifdefined\shortversion {\color{blue}
\begin{proof}
See~\cite{reportwithproof} for a detailed proof.
\end{proof} }
\fi
\ifdefined\longversion
\begin{proof} First, note that
\begin{align*}
\frac{\partial}{\partial \tau_\ell} &\int_0^T\exp(\vartheta(C_zx(t;(\tau_i)_{i=1}^m,((\alpha_{i,k})_{i=1}^{n_u})_{k=1}^K)-d_z))\mathrm{d}t\\
=&\int_0^T\vartheta\exp(\vartheta(C_zx(t;(\tau_i)_{i=1}^m,((\alpha_{i,k})_{i=1}^{n_u})_{k=1}^K)-d_z))\\
&\hspace{.2in}\times C_z \bigg[\frac{\partial}{\partial \tau_\ell}x(t;(\tau_i)_{i=1}^m,((\alpha_{i,k})_{i=1}^{n_u})_{k=1}^K)\bigg]\mathrm{d}t.
\end{align*}
Similarly, we have
\begin{align*}
\frac{\partial}{\partial \alpha_{j,\ell}} &\int_0^T\exp(\vartheta(C_zx(t;(\tau_i)_{i=1}^m,((\alpha_{i,k})_{i=1}^{n_u})_{k=1}^K)-d_z))\mathrm{d}t\\
=&\int_0^T\vartheta\exp(\vartheta(C_zx(t;(\tau_i)_{i=1}^m,((\alpha_{i,k})_{i=1}^{n_u})_{k=1}^K)-d_z))\\
&\hspace{.2in}\times C_z \bigg[\frac{\partial}{\partial \alpha_{j,\ell}}x(t;(\tau_i)_{i=1}^m,((\alpha_{i,k})_{i=1}^{n_u})_{k=1}^K)\bigg]\mathrm{d}t.
\end{align*}
The rest of the proof follows from simple algebraic manipulations.
\end{proof}
\fi

Algorithm~\ref{alg:1} may be used with the gradients in Lemma~\ref{lemma:2} to find a local solution of the optimization problem in~\eqref{eqn:optim:7}. Moreover if, after finding the optimal solution for a given $\vartheta$, the state constraints were violated at an intolerable level, we may sequentially increase $\vartheta$ and solve the problem until we get acceptable performance.

\begin{table}
\caption{\label{table:1} Numerical parameters used in the simulation.} 
\begin{tabular}{|c|c|c|c|c|c|c|}
\hline
& $c_{\mathrm{in},i}$ & $c_{\mathrm{out},i}$ & $t_{\mathrm{d},i}$ & $\kappa_i$ & $\phi_i$ & $\rho_i$ \\ \hline
$i=1$ & $0.0546$ & $0.0363$ & $5$ & $0.0103$ & $71.820$ & $8.510$ \\ \hline
$i=2$ & $0.0173$ & $0.0258$ & $6$ & $0.0084$ & $141.27$ & $16.74$ \\ \hline
\end{tabular}
\vspace{-.2in}
\end{table}

\begin{figure}\centering
\includegraphics[width=0.9\linewidth]{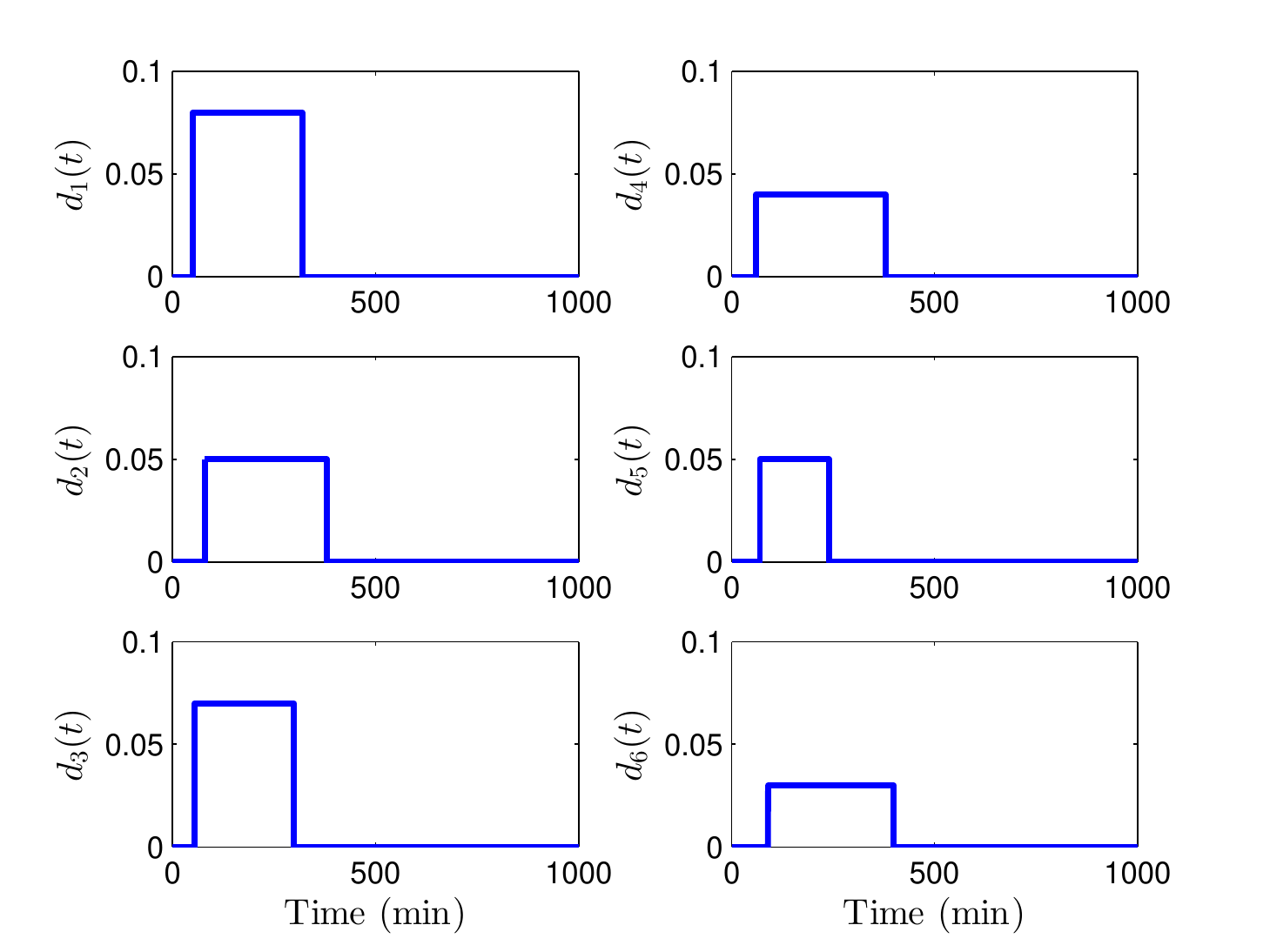}
\caption{\label{fig:demand} The demand by the farms as requested (without the scheduling delays).}
\vspace{-.3in}
\end{figure}

\begin{figure}
\vspace{-.22in}
\centering
\includegraphics[width=1.0\linewidth]{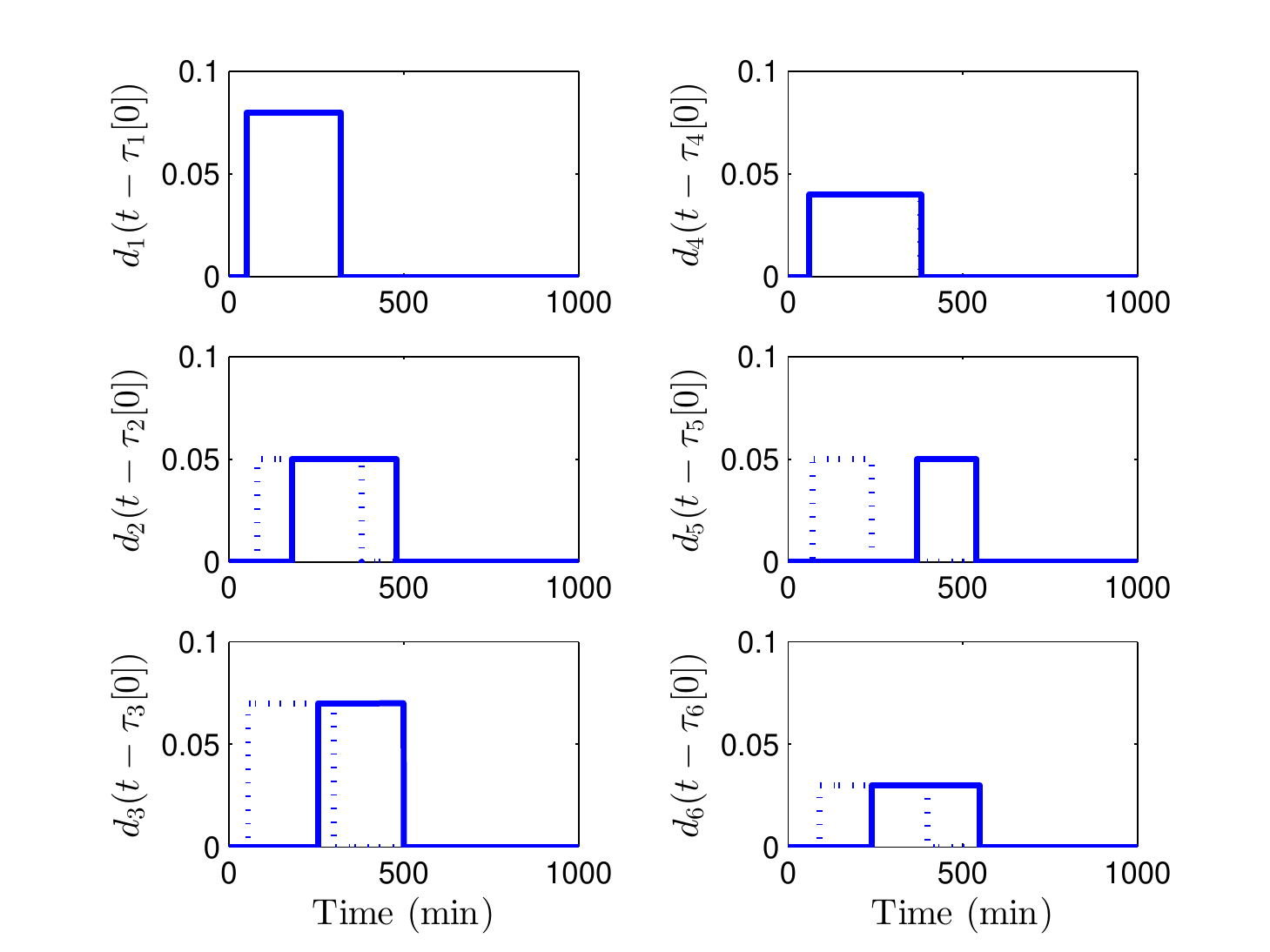}
\caption{\label{fig:initial} The shifted demands at the initialization of the algorithm. The dotted curve demonstrates the requests and the solid curve demonstrates their shifted counterpart.  }
\vspace{-.25in}
\end{figure}
\section{Numerical Example} \label{sec:numerical}
\vspace{-.04in}
In this section, we illustrate the applicability of the algorithms on a water channel with two pools. The numerical example is borrowed from~\cite{5400193}. Each pool is modelled as\vspace{-.04in}
\begin{align*}
y_i(s)=\frac{c_{\mathrm{in},i}}{s} e^{-t_{\mathrm{d},i}s}q_i(s)
-\frac{c_{\mathrm{out},i}}{s}q_{i+1}(s)-\frac{c_{\mathrm{out},i}}{s}\zeta_i(s),
\end{align*}
where $c_{\mathrm{in},i}$ and $c_{\mathrm{out},i}$ are discharge rates determined by the physical characteristics of the gates used to set the flow between neighbouring pools, and $t_{\mathrm{d},i}$ is the delay associated with the transport of water along the pool. Here, $\zeta_i(s)$ denotes the overall off-take flow load on pool $i$, that is, all the water supplied to the farms connected to this pool. Moreover, $q_i(s)$ is the flow of water from pool $i-1$ to pool $i$ and $y_i(s)$ denotes the water level in pool $i$. For the purpose of this example, we replace the delays with their first-order Pad\'{e} approximation\footnote{Note that the choice of a first-order Pad\'{e} approximation is justifiable as the pool delays are all parts of closed-loops (with local controllers), with loop-gain cross-overs that are sufficiently small to make the overall closed-loop behaviour insensitive to the approximation error~\cite{MichaelCSM}.}. Each pool is controlled, locally, by
$$
q_i(s)=\frac{\kappa_i(\phi_i s+1)}{s(\rho_i s+1)} (u_i(s)-y_i(s)),
$$
where $\kappa_i$, $\phi_i$, and $\rho_i$ are appropriately selected control parameters. Furthermore, $u_i(s)$ denotes the water-level reference signal of pool $i$. Table~\ref{table:1} shows the parameters used in this example. 
The state constraints are as follows
$9.4 \leq y_1(t)\leq 9.7$ and $9.5 \leq y_2(t)\leq 9.7$.
Finally, throughout this example, we fix 
$
u_0=[9.50,9.55]^\top.
$

\begin{figure}\centering
\includegraphics[width=1.0\linewidth]{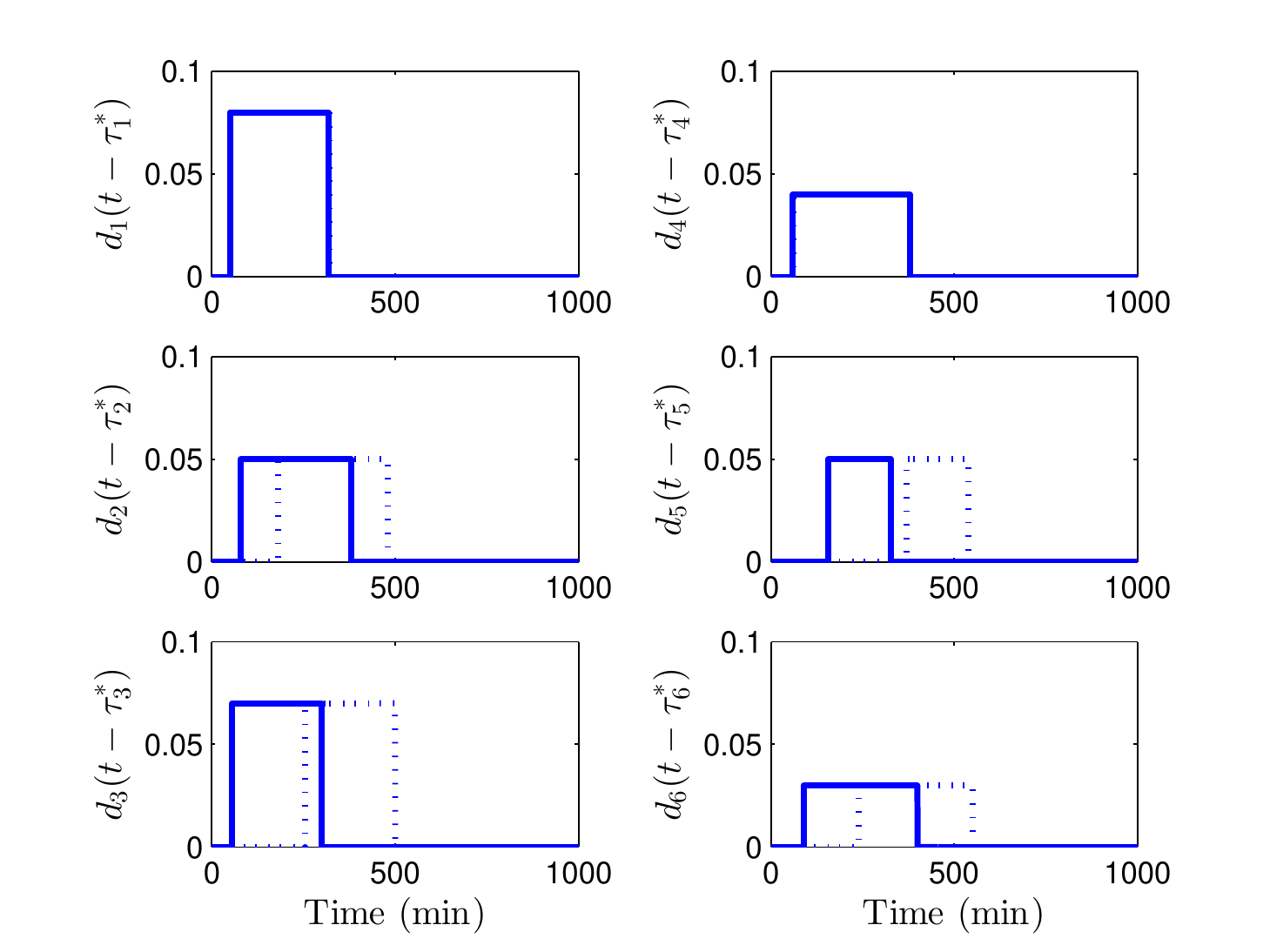}
\caption{\label{fig:ZOH_delay} The shifted demands for the local solution recovered by Algorithm~\ref{alg:1} with $\epsilon=0.1$. The dotted curve demonstrates the shifted demands at the initialization and the solid curve demonstrates the shifted demands at the locally optimal solution. }
\vspace{-.2in}
\end{figure}

\begin{figure}\centering
\includegraphics[width=0.9\linewidth]{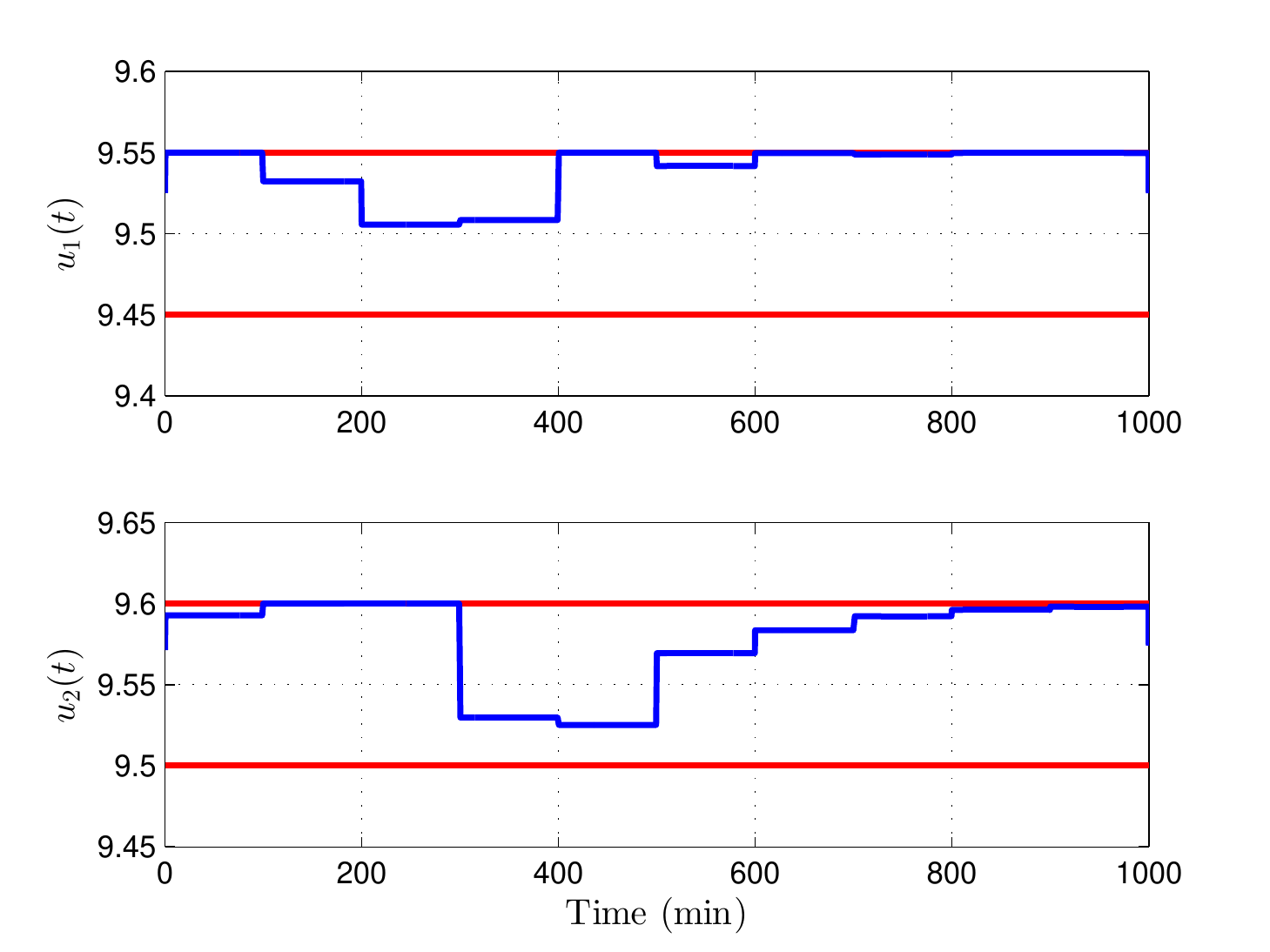}
\caption{\label{fig:ZOH_control} The reference signal for the local solution recovered by Algorithm~\ref{alg:1} with $\epsilon=0.1$. The red lines show the boundary of the feasible region. }
\vspace{-.2in}
\end{figure}

Figure~\ref{fig:demand} illustrates the requested demands of the farms. Here, $(v_i(t))_{i=1}^3$ and $(v_i(t))_{i=4}^6$, respectively, denote demands for pool $1$ and $2$. Let us select linear penalty functions $h_i(\tau_i)=\tau_i$ for all $i$. Moreover, assume that the reference signal should belong to a bounded region captured by $$u_0-\begin{bmatrix}
0.05 \\ 0.05
\end{bmatrix} \leq u(t)\leq u_0+\begin{bmatrix}
0.05 \\ 0.05
\end{bmatrix}.$$ 
Note that without these control input constraints, one can schedule all the loads without any delay but with large control input deviations. In this example, we select $\underline{\tau}_i=0$, $\forall i$, which means that we can only shift demands forward. Let us also select $\overline{\tau}_i=300\,\mathrm{min}$ for all $i$.

\begin{figure}\centering
\includegraphics[width=0.9\linewidth]{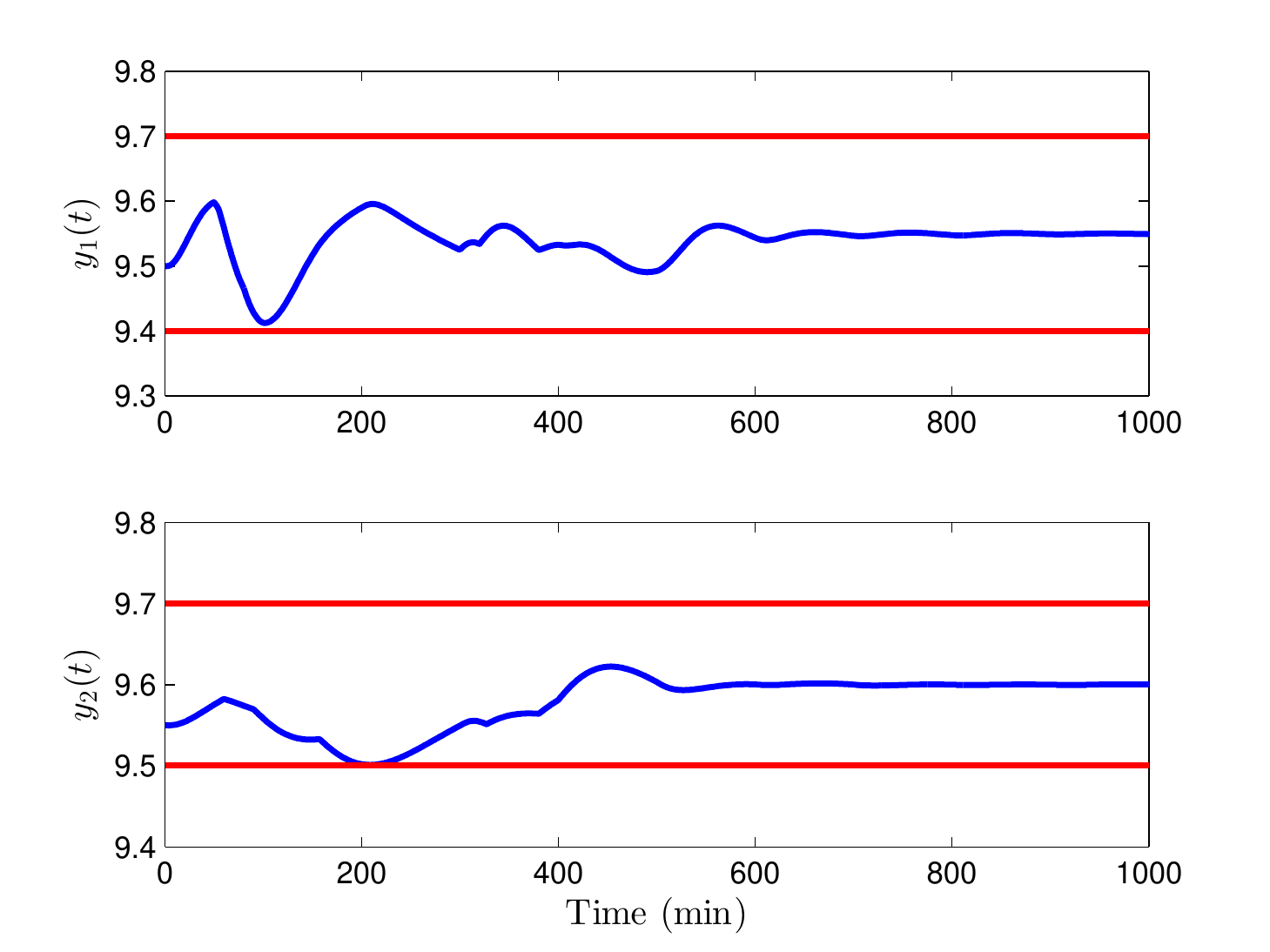}
\caption{\label{fig:ZOH_state} The outputs for the local solution recovered by Algorithm~\ref{alg:1} with $\epsilon=0.1$. The red lines show the boundary of the feasible region. }
\vspace{-.2in}
\end{figure}

\begin{figure}\centering
\includegraphics[width=1.0\linewidth]{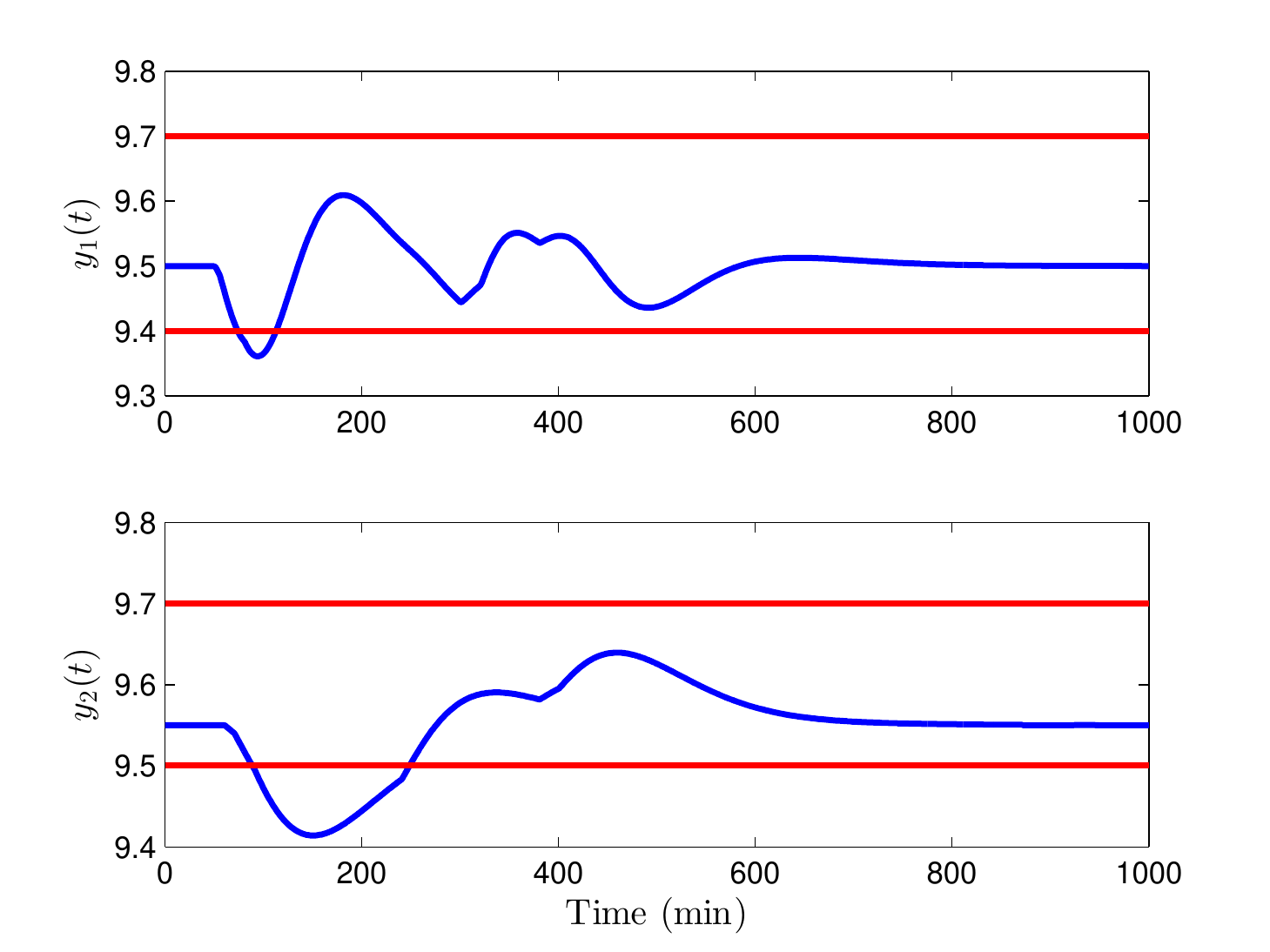}
\caption{\label{fig:Soft_intial} The output of the system when all the decision variables (the control inputs and scheduling delays) are set equal to zero. }
\vspace{-.2in}
\end{figure}

\begin{figure}\centering
\includegraphics[width=1.0\linewidth]{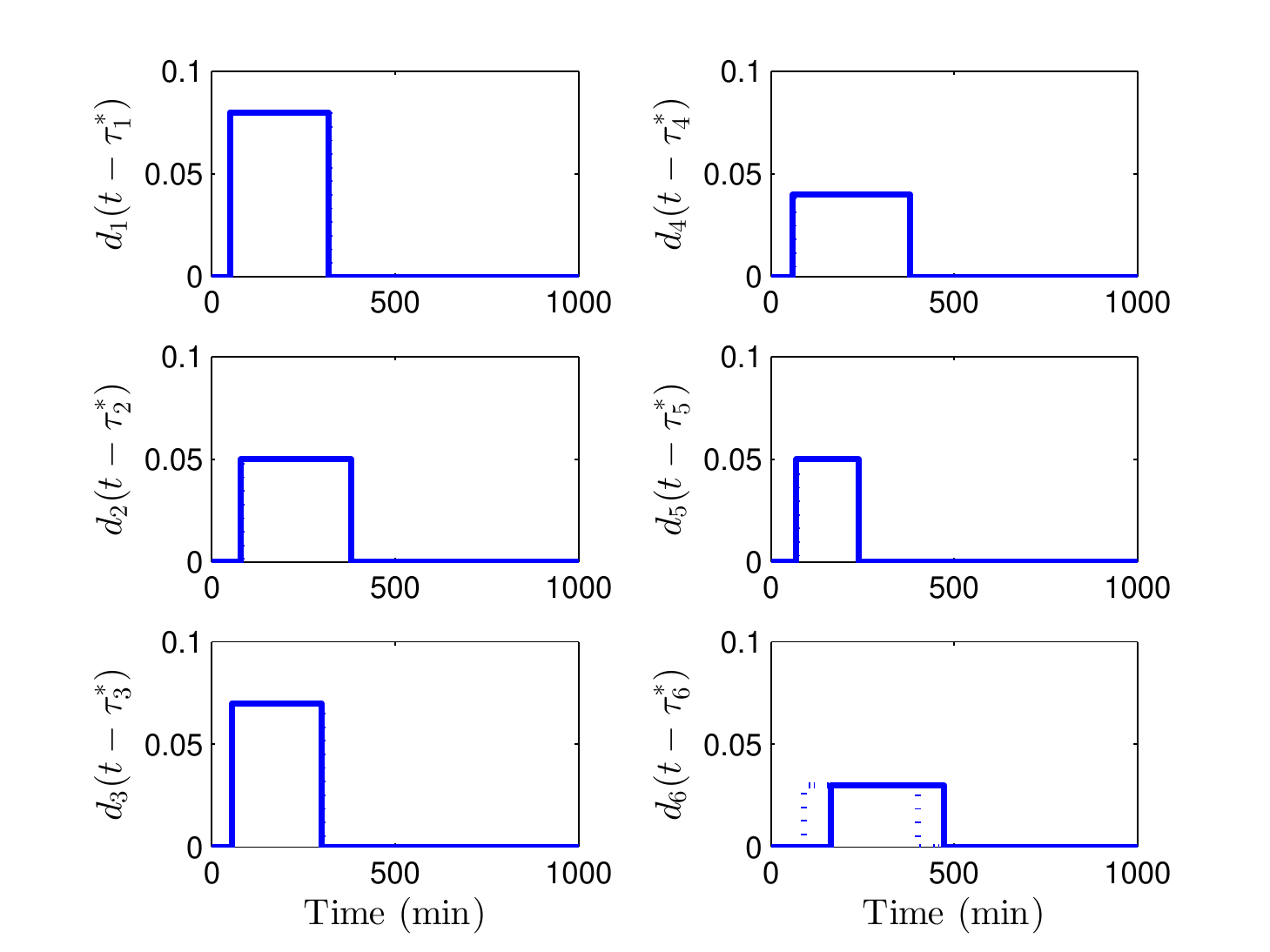}
\caption{\label{fig:Soft_delay} The shifted demands for the local solution recovered by the proposed algorithm in Section~\ref{sec:soft} with $\vartheta=100$. The dotted curve demonstrates the shifted demands at the initialization and the solid curve demonstrates the shifted demands at the suboptimal solution. }
\vspace{-.2in}
\end{figure}

First, we use Algorithm~\ref{alg:1} to extract a reasonable schedule by shifting these demands. This algorithm requires a feasible starting point, which can be constructed by shifting the demands (to be somewhat distant from each other). Figure~\ref{fig:initial} illustrates the shifted demands at this initialization (solid curve) as well as the original requests (dotted curves) for comparison. Let us fix $\epsilon=0.1$. Figure~\ref{fig:ZOH_delay} shows the shifted demands for the local solution recovered by the proposed algorithm. As we expect, the delays are significantly smaller in comparison to the initialization. Figure~\ref{fig:ZOH_control} portrays the reference signals for the solution and Figure~\ref{fig:ZOH_state} illustrates the outputs. Evidently, the output stays in the desired region. Although powerful, Algorithm~\ref{alg:1} requires a feasible initial condition that may not be easy to find. Therefore, in the rest of this section, we study the method presented in Section~\ref{sec:soft}.

\begin{figure}\centering
\includegraphics[width=0.9\linewidth]{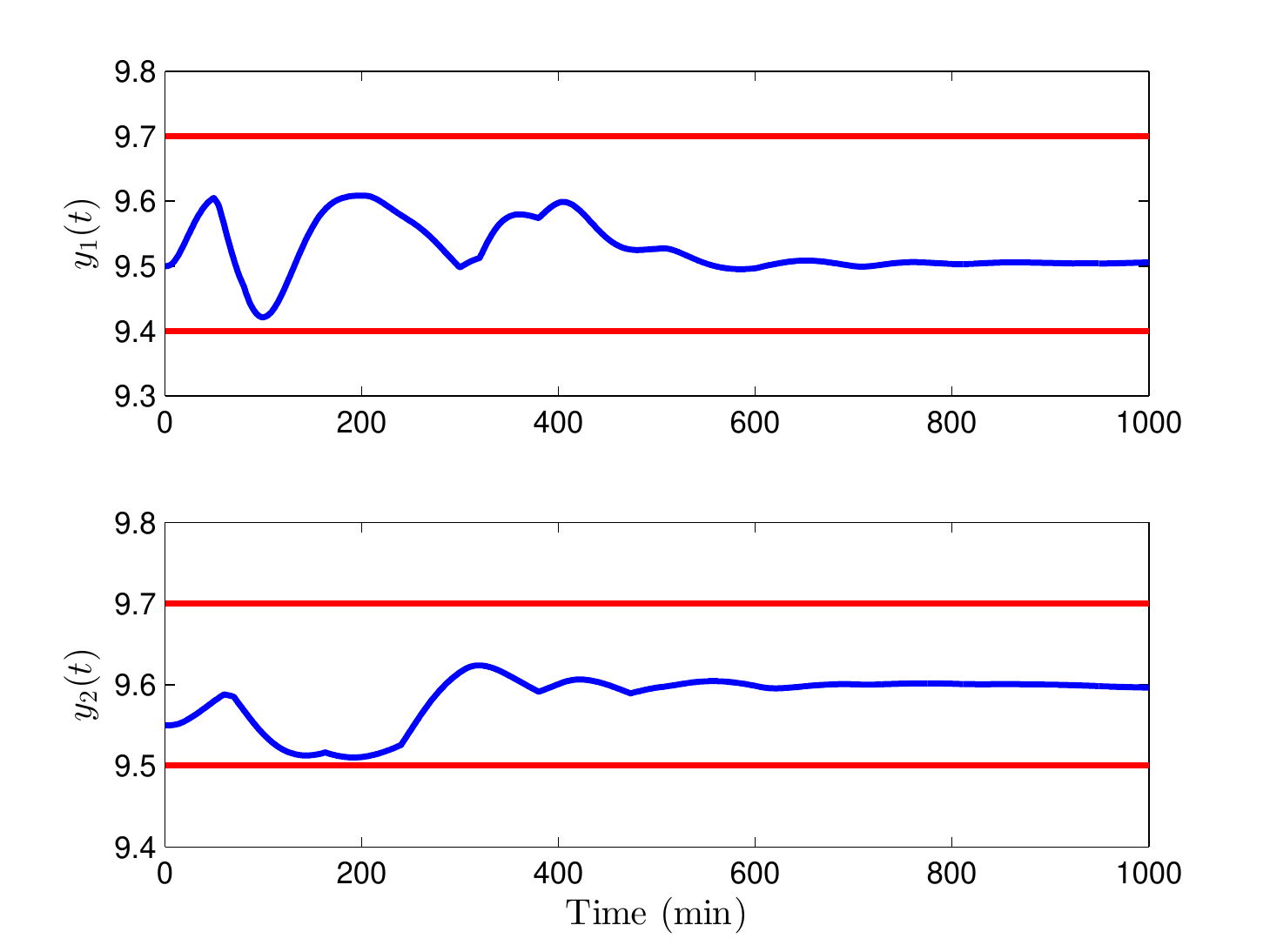}
\caption{\label{fig:Soft_state} The outputs for the locally optimal solution recovered by the proposed algorithm in Section~\ref{sec:soft} with $\vartheta=100$. The red lines show the boundary of the feasible region. }
\vspace{-.2in}
\end{figure}

\begin{figure}\centering
\includegraphics[width=0.9\linewidth]{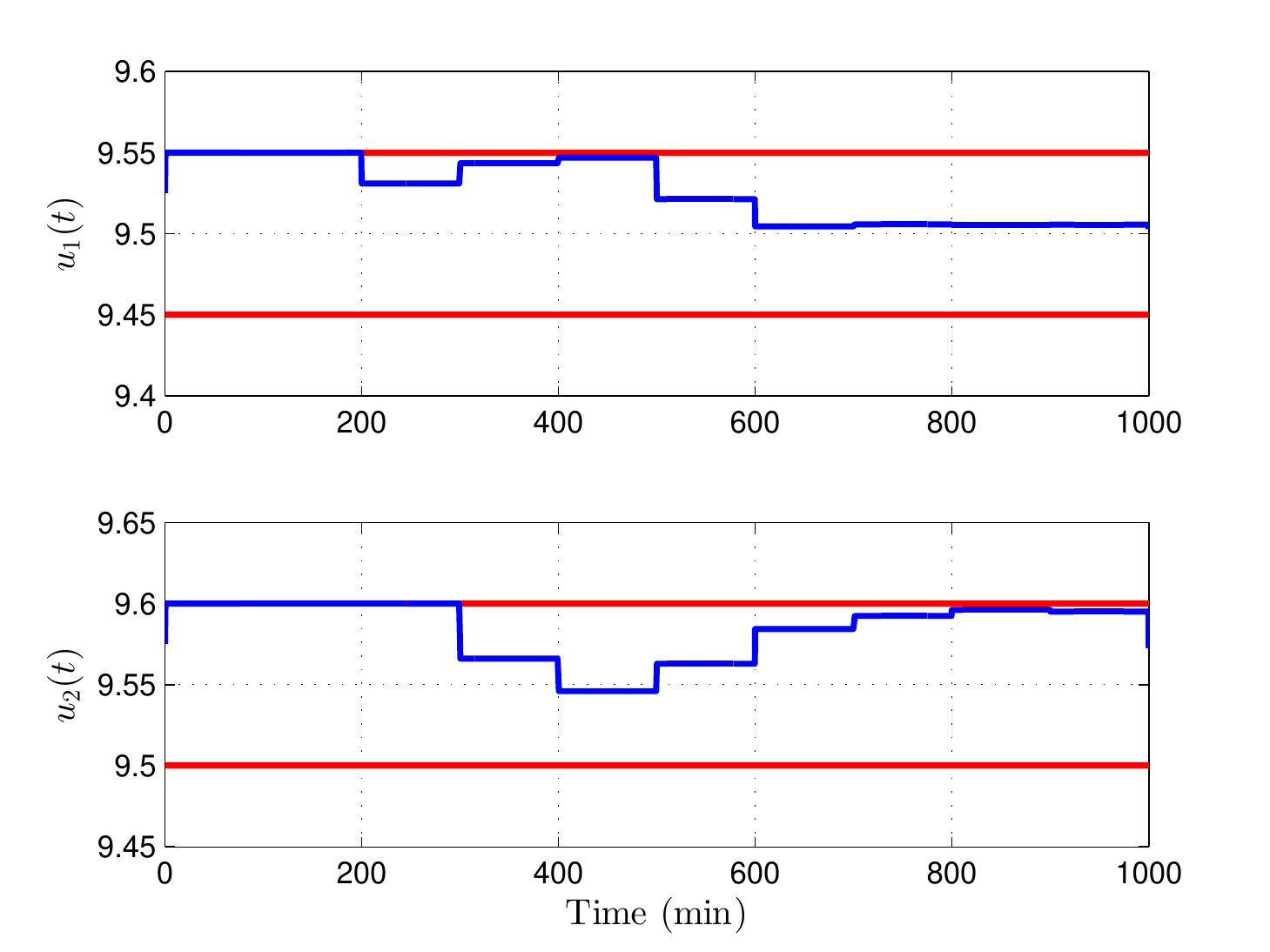}
\caption{\label{fig:Soft_control} The reference signal for the locally optimal solution recovered by the proposed algorithm in Section~\ref{sec:soft} with $\vartheta=100$. The red lines show the boundary of the feasible region. }
\vspace{-.2in}
\end{figure}

\begin{figure}\centering
\includegraphics[width=0.9\linewidth]{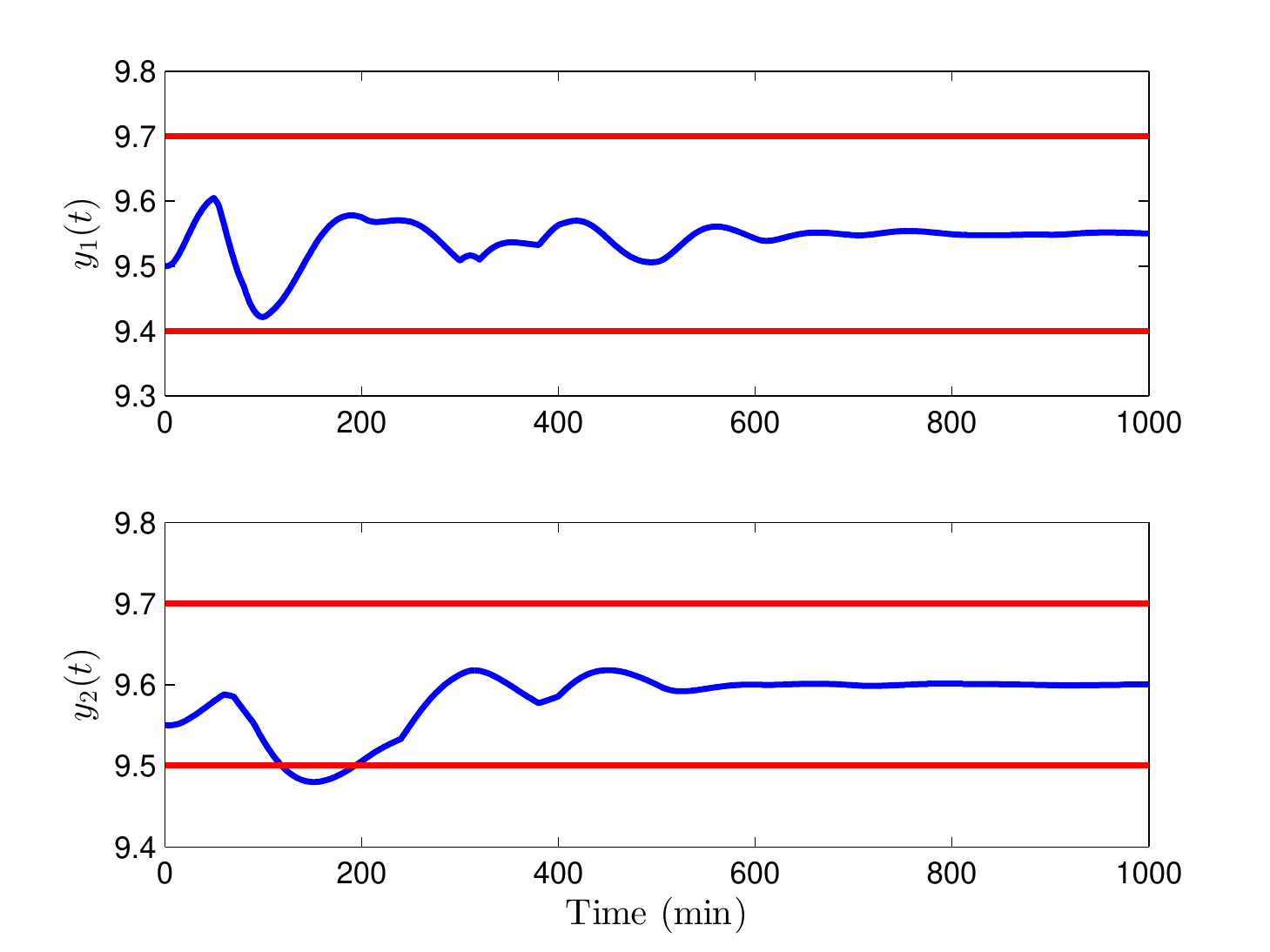}
\caption{\label{fig:Soft_state_violate} The outputs for the locally optimal solution recovered by the proposed algorithm in Section~\ref{sec:soft} with $\vartheta=10$. The red lines show the boundary of the feasible region. }
\vspace{-.2in}
\end{figure}

Let us select $\vartheta=100$. Figure~\ref{fig:Soft_delay} illustrates the shifted demands for the locally optimal solution recovered by the proposed algorithm. At the starting point (fed to the algorithm), all the decision variables are selected to be equal to zero for which the state of the system does not stay in the feasible region; see Figure~\ref{fig:Soft_intial}. Figures~\ref{fig:Soft_state} and~\ref{fig:Soft_control} show the output and the control  for the local solution recovered by the proposed algorithm in Section~\ref{sec:soft}. Interestingly, all the constraints are satisfied, which is because of the large value of $\vartheta$. If we reduce $\vartheta$ to be equal to $10$, the outputs violate the constraints on the state as shown in Figure~\ref{fig:Soft_state_violate}. Evidently, by comparing Figures~\ref{fig:Soft_intial}, \ref{fig:Soft_state}, and~\ref{fig:Soft_state_violate}, we can see that by increasing $\vartheta$, the constraint violations are becoming more infrequent (until they do not occur at all).

\section{Conclusions} \label{sec:conclusions}
In this paper, we presented numerical algorithms for scheduling demands on continuous-time linear time-invariant systems. The rigidity of the demands dictated that we can only shift them  back-and-forth in time (and cannot change their shapes). The first algorithm used log-barrier functions to include the state constraints in the cost function. The second algorithm considered the state constraints as soft constraints and added a penalty function for the constraint violations to the cost function. Future research can focus on constructing a market mechanism for achieving the optimal schedule based on the customers preferences. We can also compute an optimality gap via finding a lower-bound for the solution of the dual of the problem corresponding to~\eqref{eqn:optim:1}.

\bibliographystyle{IEEEtran}
\bibliography{citation}

\end{document}